\documentclass[11pt]{amsart}
\usepackage{amsmath,amssymb}
\usepackage{geometry}
\usepackage{amsmath, amsfonts,amsthm,amssymb,enumerate,amscd}

\usepackage[active]{srcltx}
\usepackage[utf8]{inputenc}
\usepackage{algorithm}
\usepackage{stmaryrd}
\usepackage{array}
\usepackage{algorithmicx}
\usepackage{algpseudocode}
\usepackage{graphicx}
\usepackage{letltxmacro}

\usepackage{color}
\usepackage[table,xcdraw,svgnames,dvipsnames]{xcolor}
\definecolor{light-salmon}{RGB}{255,140,120}
\usepackage[colorlinks]{hyperref}
\hypersetup{
	allbordercolors=.,
	citecolor=DarkGreen,
	linkcolor=DarkRed,
	urlcolor=NavyBlue
}

\graphicspath{{./pics/}}

\usepackage{enumitem}

\numberwithin{equation}{section}
\theoremstyle{plain}
\newtheorem{thm}{Theorem}

\newtheorem{cor}[thm]{Corollary}
\newtheorem{prop}[thm]{Proposition}

\theoremstyle{definition}
\newtheorem{rem}[thm]{Remark}

\renewcommand{\Bbb}{\mathbb}
\newcommand{\bo}[1]{{\bf #1}}

\makeatletter
\DeclareFontFamily{U}{tipa}{}
\DeclareFontShape{U}{tipa}{m}{n}{<->tipa10}{}
\newcommand{\arc@char}{{\usefont{U}{tipa}{m}{n}\symbol{62}}}%

\newcommand{\arc}[1]{\mathpalette\arc@arc{#1}}

\newcommand{\arc@arc}[2]{%
  \sbox0{$\m@th#1#2$}%
  \vbox{
    \hbox{\resizebox{\wd0}{\height}{\arc@char}}
    \nointerlineskip
    \box0
  }%
}
\makeatother

\title[On the Blaschke-Lebesgue theorem for the Cheeger constant]{On the Blaschke-Lebesgue theorem for the Cheeger constant via areas and perimeters of inner parallel sets}

\author{Beniamin Bogosel}

\begin{document}
	
\maketitle

\begin{abstract}
	The first main result presented in the paper shows that the perimeters of inner parallel sets of planar shapes having a given constant width are minimal for the Reuleaux triangles. This implies that the areas of inner parallel sets and, consequently, the inverse of the Cheeger constant are also minimal for the Reuleaux triangles. Proofs use elementary geometry arguments and are based on direct comparisons between general constant width shapes and the Reuleaux triangle.
\end{abstract}

\bo{Keywords:}  constant width, inner parallel sets, Cheeger constant


\section{Introduction}

Considering a convex domain $\Omega$ in the plane, i.e. a convex and closed set, a supporting line is a line which intersects $\Omega$ but does not separate any two points in $\Omega$. For a smooth region of $\partial \Omega$, supporting lines coincide with tangent lines. Given an orientation $\theta\in [0,2\pi]$ in the plane, the distance between the two supporting lines orthogonal to $\theta$ is called the width $w(\Omega,\theta)$ of $\Omega$ in direction $\theta$. Shapes which have the same width for any direction $\theta \in [0,2\pi]$ are called \emph{constant-width} shapes. The most basic examples of constant width shapes are the disk and the Reuleaux triangle, the intersection of three disks having radius one with centers at the vertices of an equilateral triangle of side length $1$. It turns out that the Reuleaux triangle is extremal for various geometric quantities and in this paper additional such results are proved. The following result attributed to Lebesgue \cite{Lebesgue} and Blaschke \cite{Blaschke1915} has attracted a lot of attention.

\bo{Blaschke-Lebesgue Theorem.} The Reuleaux triangle minimizes the area among shapes having a given constant width.

This result has many different proofs besides the original ones. Various geometric arguments are given in the book of Yaglom and Boltyanskii \cite{yaglom-boltjanskii} which has a full chapter on constant width shapes. Chakerian gives a surprising proof in \cite{Chakerian} using circumscribed hexagons. Ghandehari uses control theory in \cite{Ghandelhari} and Harrell uses variational techniques in \cite{Harrell}.

It is known that the Reuleaux triangle also minimizes the inradius and maximizes the circumradius \cite[Chapter 7]{yaglom-boltjanskii}. More recently, other optimization problems were studied in the class of shapes of constant width. The minimization of the eigenvalues of the Dirichlet-Laplace operator under diameter constraint yields optimal shapes which have constant width \cite{BHL18}. The maximization of these eigenvalues is also well posed in the class of constant width shapes and numerical simulations show that it is likely that the Reuleaux triangle is again the optimal shape \cite{Bogosel-Convex}.

Colesanti proves in \cite{Colesanti} that the Brunn-Minkowski inequality holds for various functionals, like the first eigenvalue of the Dirichlet $p$-Laplace operator $\lambda_{1,p}$. In particular, since the Cheeger constant $h(\Omega)=\inf\{ |\partial X|/|X| \text{ such that } X \subset \Omega\}$, is obtained as the limit of the eigenvalues $\lambda_{1,p}(\Omega)$, Brunn-Minkowski inequality also extends to this case. Basic convexity arguments show that shapes of constant width maximizing $\lambda_{1,p}$ or the Cheeger constant must be indecomposable, i.e. they cannot be written as the Minkowski sum of two non-homothetic constant width shapes. The description of such indecomposable bodies goes beyond the scope of this article, but we may cite the example of Reuleaux polygons. Similar techniques involving the minimization of concave functionals are illustrated in \cite{bucur-concave}. Unfortunately, in most cases underlined above, indecomposability is not enough to conclude that the Reuleaux triangle is optimal. 

The first eigenvalue of the $L^\infty$ Laplacian of $\Omega$ is $1/r(\Omega)$, the inverse of the inradius, and is maximized by the Reuleaux triangle $R$. In view of the numerical results shown in \cite{Bogosel-Convex} concerning the $L^2$-Laplacian, it is therefore natural to conjecture that the other extremal case, the Cheeger constant, corresponding to the $1$-Laplacian is also maximized by $R$. This was proved recently by Henrot and Lucardesi in \cite{henrot-lucardesi} using techniques from shape optimization, namely optimality conditions related to the shape derivative, which are verified by the minimizer.

In this paper a different approach which is completely geometric in nature gives a proof of the same result. The Cheeger constant of planar convex sets has a characterization based on the area of the inner parallel sets given by Lachand-Robert and Kawohl in \cite{kawohl-lachand-robert}. More precisely, if 
\begin{equation}
\Omega_{-t}= \{ x \in \Omega: d(x,\partial \Omega) \geq t\}
\label{eq:def-inner-parallel}
\end{equation}
is the inner parallel set at distance $t$ from the boundary of the convex domain $\Omega$ then $h(\Omega) = 1/t$, where $t$ verifies $|\Omega_{-t}|=\pi t^2$. Contrary to outer parallel sets of a convex set where Steiner's formula provides a polynomial expression for the area in terms of the distance to the boundary (see for example \cite[Chapter 4]{schneider}), no such formula exists for inner parallel sets. Details for the polygonal case are given in \cite{kawohl-lachand-robert} and in this work, similar computations are made for Reuleaux polygons.

In this article new proofs are given for the following results:
\begin{enumerate}[label=(\roman*),noitemsep,topsep=0pt]
	\item Given $t\in (0,1-\sqrt{3}/{3}]$, the Reuleaux triangle minimizes the \bo{perimeter} of $\Omega_{-t}$ among shapes with constant width equal to one. The result is proved in Theorem \ref{thm:perim-parallel}. This result is new up to the author's knowledge.
	\item Given $t\in [0,1-\sqrt{3}/{3}]$, the Reuleaux triangle minimizes the \bo{area} of $\Omega_{-t}$ among shapes with constant width equal to one. The result is proved in Theorem \ref{thm:perim-parallel}. Moreover, this provides another proof of the Blaschke-Lebesgue theorem, as a direct consequence of Theorem \ref{thm:perim-parallel} and the minimality of the inradius. The article \cite{disk-polygons-BL} extends the Blaschke-Lebesgue theorem to disk-polygons (intersection of equal disks in the plane) and its results, although not explicitly stated, imply Theorem \ref{thm:area-parallel}. 
	\item The Reuleaux triangle maximizes the \bo{Cheeger constant} among shapes with constant width equal to one. The result is proved in Theorem \ref{thm:cheeger} and is a direct corollary of the characterization given in \cite{kawohl-lachand-robert} and Theorem \ref{thm:area-parallel}. The result was initially proved in \cite{henrot-lucardesi} using optimality conditions related to the shape derivative. The proof following from the results of this paper is new and elementary.
\end{enumerate}

The paper is organized as follows. Section \ref{sec:prelim} recalls basic aspects regarding constant width shapes, disk-polygons and gives a geometric lemma regarding the convexity of a function used in the sequel. Section \ref{sec:area-perim} presents proofs of the Blaschke-Lebesgue theorem for perimeters and areas of inner parallel sets. Section \ref{sec:cheeger} deals with the proof of the maximality of the Reuleaux triangle for the Cheeger constant.

\section{Preliminaries}
\label{sec:prelim}
\subsection{Planar constant width shapes}
A two dimensional convex set $K$ has constant width if the distance between any two parallel supporting lines is constant regardless of the orientation of these lines. Recall that a supporting line intersects the boundary of $K$ and leaves the shape $K$ in one of the half-planes generated by this line. The circle is the most obvious example of a shape having constant width. However, many more such shapes exist. In the rest of the article the width of the shapes considered is always equal to $1$.

The most famous examples of constant width shapes are the Reuleaux triangle and more generally, Reuleaux polygons. These shapes are not just mathematical curiosities, but have various applications \cite[Chapter 18]{gardner}. Reuleaux triangles are used in rotary engine design and square hole drilling machines, while the twenty pence British coin is a Reuleaux heptagon.

The Reuleaux triangle is defined as the intersection of the disks of radius $1$ centered at the vertices of an equilateral triangle of edge length $1$. Reuleaux polygons can be constructed using a similar procedure as shown in \cite[Section 8.1]{bodies_of_constant_width}, for example. It is classical that any constant width shape can be approximated arbitrarily well with a Reuleaux polygon \cite[Theorem 8.1.1]{bodies_of_constant_width}, \cite[Chapter 7]{yaglom-boltjanskii}.

The book \emph{Convex Figures} by Yaglom and Boltyanskii \cite{yaglom-boltjanskii} has a whole chapter dedicated to such shapes. Results are presented in the form of exercises with proofs using elementary geometry aspects. Let us recall a few of these results, relevant to this work.

\begin{prop} A planar shape with constant width $1$ has the following properties:
\begin{enumerate}[label=(\roman*),noitemsep,topsep=0pt]
\item Any shape of constant width can be approximated arbitrarily well by Reuleaux polygons having the same width, with respect to the Hausdorff distance (\cite[Theorem 8.1.1]{bodies_of_constant_width}).
\item The interior angle at a corner point in a constant width shape cannot be less than $2\pi/3$. Moreover, if a constant width curve has an interior angle equal to $2\pi/3$ then this curve is a Reuleaux triangle. \cite[Exercise 7-9]{yaglom-boltjanskii}  The presence of a vertex $A$ corresponding to a corner point of angle $\beta$ in the boundary of $K$ implies the existence of a circular arc of radius opposite to $A$ in $\partial K$ of angle $\pi-\beta$. Conversely, any circular arc of radius $1$ in the boundary corresponds to a corner point. In particular, any circular arc of radius $1$ in the boundary of $K$ has length or angle measure at most $\pi/3$.

\item The inscribed and circumscribed disks to a constant width shape are concentric and the sum of their rays is equal to the constant width. \cite[Exercise 7-13]{yaglom-boltjanskii} 

\item The Reuleaux triangle is the curve of constant width with the greatest circumradius, therefore having the smallest inradius. \cite[Exercise 7-14]{yaglom-boltjanskii}  

\item The perimeter of curves having constant width $1$ is equal to $\pi$. In particular, since any Reuleaux polygon has boundary consisting of a series of arcs of circles of radii equal to $1$, the sum of the subtended angles of these arcs is equal to $\pi$.

\item The area of a constant width shape is minimized by the Reuleaux triangle (Blaschke-Lebesgue theorem) and is maximized by the disk. Instructive proofs of the last two points are given in the quoted book, based on circumscribed equiangular polygons. \cite[Exercise 7-12]{yaglom-boltjanskii}
\end{enumerate}
\label{prop:properties}
\end{prop}

The Blaschke-Lebesgue theorem asserts that the Reuleaux polygon minimizes the area among shapes of fixed constant width. The result is attributed to Blaschke \cite{Blaschke1915} and Lebesgue \cite{Lebesgue}. Many other proofs strategies have been used to prove the result, among which we mention \cite{Ghandelhari} using an optimal control formulation and \cite{Harrell} using variational techniques. An overview of the existing bibliography is given in \cite[Theorem 12.1.5]{bodies_of_constant_width}. The results of this paper give yet another different proof of this result.

In general, Reuleaux polygons are assumed to have an odd number of arcs. A careful analysis shows that if an even number of arcs are present, then some consecutive arcs correspond to the same center and can be merged. A Reuleaux polygon is \emph{regular} if the centers of the arcs forming its boundary are the vertices of a regular polygon. In \cite{Firey1960} Firey shows that among Reuleaux polygons with fixed number of arcs, the regular one has the largest area. The proof, based on estimates related to areas of parallel inner sets and convexity arguments, partially inspired some of the methods used in this paper.

Given a Reuleaux polygon with $n=2k+1$ vertices, its boundary is made of arcs of circle subtending angles $\theta_i$, $i=0,...,n-1$ which verify $\sum_{i=0}^{n-1} \theta_i = \pi$. Following Proposition \ref{prop:properties}-(ii) we also have $\theta_i \in [0,\pi/3]$. 

\subsection{Disk Polygons}
\label{sec:disk-poly}

Reuleaux polygons have the particularity that they are convex sets whose boundaries are made of finitely many arcs of circles of having the same radius. More precisely, Reuleaux polygons are intersection of disks having the same radius. For simplicity, in the following, we assume that disks have radius one unless stated otherwise. 

For dealing with inner parallel sets of Reuleaux polygons it is useful to define an even more general class of convex sets, namely the \emph{disk polygons}, as the intersection of a finite number of disks of radius one. This concept is natural and was introduced previously, for example in \cite{disk-polygons}. Moreover, in \cite{disk-polygons-BL} it was proved that the Reuleaux triangle (which is obviously a disk polygon) minimizes the area among all disk polygons whose centers are at distance at most $1$ apart. For the purpose of this article, we only need to investigate disk-polygons which already contain a Reuleaux polygon of width $1$. Basic properties of disk polygons are recalled below.

Consider $\Omega = \bigcap_{i=0}^{N-1} D_i$ a disk polygon, where $D_i$ are disks of radius one. Suppose also that every disk \emph{contributes} to $\Omega$, i.e. the family $D_i$ is minimal, in particular, no disk is duplicated. Denote by $\Gamma_i$, $i=0,...,N-1$ the arcs defining the boundary of $D$, having lengths $\theta_i$, respectively. Again, the arcs $\Gamma_i,\Gamma_j$ are assumed to belong to different disks if $i \neq j$. The arcs $\Gamma_i$ have the extremities at vertices $v_i, v_{i+1}$, where indices are taken modulo $n$. At each vertex $v_i$, the arcs $\Gamma_i, \Gamma_{i-1}$ meet with the \emph{turning angle} $\beta_i$ (the angle made at $v_i$ by the tangent vectors at $\Gamma_{i-1},\Gamma_i$  in the trigonometric sense). A couple of properties of interest for the sequel of the paper are presented below. 

\begin{prop}
	\begin{enumerate}[label=(\roman*),noitemsep,topsep=0pt]
		\item The sum of the lengths $\theta_i$ of the arcs $\Gamma_i$ and of the turning angles $\beta_i$ is equal to $2\pi$: 
		\[ \sum_{i=0}^{N-1} (\theta_i+\beta_i)=2\pi.\]
		\item If $\Omega$ is a disk polygon which contains a shape $\Omega'$ of constant width equal to $1$ then:
		\begin{itemize}[noitemsep]
			\item The minimal width of $\Omega$ is at least $1$.
			\item The perimeter of $\Omega$ is at least $\pi$: $\sum_{i=0}^{N-1} \theta_i \geq \pi$.
			\item For every turning angle we have $\beta_i \in [0,\pi/3]$, $i=0,...,N$.
			\item The inradius $r(\Omega)$ verifies $r(\Omega) \geq r(R) = \frac{3-\sqrt{3}}{2}$, where $R$ is the Reuleaux triangle of width $1$.
		\end{itemize}
	\end{enumerate}
	\label{prop:disk-poly}
\end{prop}

\emph{Proof:} (i) It is enough to follow a turning supporting line around $\Omega$. Each arc $\Gamma_i$ turns the line with an angle equal to $\theta_i$. Each vertex contributes with the turning angle $\beta_i$. 

(ii) Any pair of parallel supporting lines to $\Omega$ generate a strip containing the constant width shape $\Omega'$. Therefore, the width of any such strip is at least equal to $1$.

The perimeter of convex sets is monotone with respect to inclusion (see for example \cite[Lemma 2.2.2]{bucur-buttazzo}). Therefore $|\partial \Omega|\geq |\partial \Omega'| = \pi$.

Consider two disks $D_i, D_{i-1}$ containing arcs $\Gamma_i, \Gamma_{i-1}$ meeting at a vertex $v_i$ having turning angle $\beta_i$. Then $\Omega \subset D_i \cap D_{i-1}$ and the minimal width of the disk intersection $D_i \cap D_{i-1}$ is at least equal to one, i.e. their centers are at distance at least $1$ apart. Therefore, the turning angle $\beta_i$ at the intersection of their boundaries is at most $\pi/3$.

The inradius is monotone with respect to inclusion so the last point follows immediately, since the Reuleaux triangle minimizes the inradius among shapes with given constant width. 
\hfill $\square$


\subsection{A geometrical lemma}
Consider a fixed $t \in (0,(3-\sqrt{3})/3)$ and a triangle $\Delta ABC$ such that $AB=1$, $BC=1-t$. Denoting with  $\angle BAC = \gamma$ denote $\alpha(\gamma) = \angle CBA$ and with $h(\gamma)$ the distance from $C$ to $AB$, the height of $\Delta ABC$ from the vertex $C$. The previous notation emphasizes that given $\gamma$, supposing the lengths of $AB, BC$ fixed and $\alpha(\gamma)$ acute, the angles $\alpha(\gamma)$ and the height $h(\gamma)$ can be expressed in terms of $\gamma$.

Equivalently, consider $A$ at the origin, $B$ with coordinates $(1,0)$ and $D(1-t)$ the circle with center $B$ and radius $1-t$. The line $\ell(\gamma)$ through $A$, making an angle $\gamma$ with $AB$ intersects $D(1-t)$ at two points, the closest one to $A$ being $C(\gamma)$. The largest value $\gamma$ for which $C(\gamma)$ is well defined corresponds to the case where $AC(\gamma)$ is tangent to $D(1-t)$, in which case, $\gamma = \arcsin(1-t)$. In the application, $t$ is bounded above by the inradius of the Reuleaux triangle, $r(R) = \frac{3-\sqrt{3}}{3}$, which implies $1-t \geq \sqrt{3}/3$. Moreover, in the application we only consider $\gamma \in [0,\pi/6]$, corresponding to half of a turning angle (see Propositions \ref{prop:properties}, \ref{prop:disk-poly}). Considering the triangle $YAB$ with $\angle YAB=\angle YBA = \pi/6$, choosing $Y$ in the first quadrant, we observe that $YA=YB=\sqrt{3}/{3}$. Therefore, $D(1-t)$ always intersects $YA$ since $ \sqrt{3}/{3} \leq 1-t \leq 1$, implying that $C(\gamma)$ is well defined, assuming $\gamma \in [0,\pi/6]$ and $t \leq (3-\sqrt{3})/3$.

\begin{figure}
	\centering 
	\includegraphics[width=0.8\textwidth]{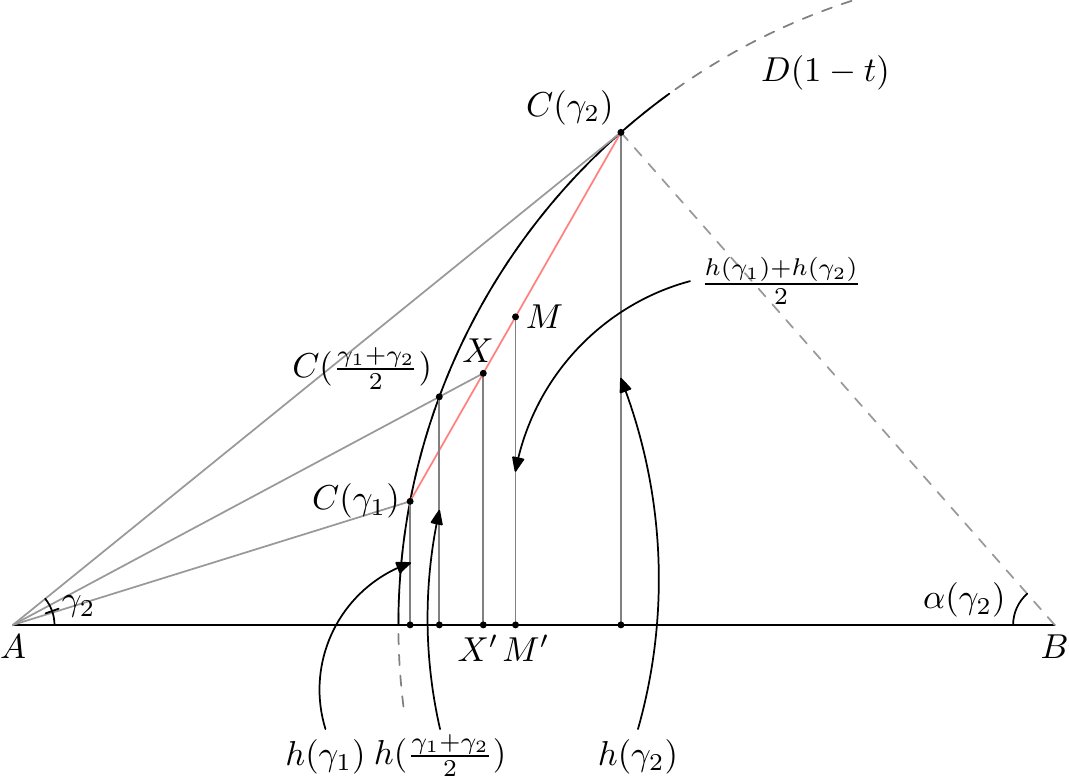}
	\caption{Configuration corresponding to Proposition \ref{prop:convexity}: geometric proof of the convexity of the application $\gamma \mapsto h(\gamma)$.}
	\label{fig:convexity}
\end{figure}

Therefore, we may define 
\[ \alpha(\gamma), h(\gamma) : [0,\pi/6]\to \Bbb{R}\]
with the properties stated above. See Figure \ref{fig:convexity} for an illustration.
The next geometrical lemma, which is fundamental for the results of the paper, shows that the dependence of $\alpha(\gamma),h(\gamma)$ in $\gamma$ is convex. 

\begin{prop}
	The applications $\gamma \mapsto \alpha(\gamma)$ and $\gamma \mapsto h(\gamma)$ are strictly increasing and strictly convex. Moreover, we also have $\alpha(\gamma)\leq \gamma$ for $\gamma \in [0,\pi/6]$.
	\label{prop:convexity}
\end{prop}

\emph{Proof:} Since $\alpha(\gamma) =\arcsin \frac{h(t)}{1-t}$ and $\arcsin$ is convex and strictly increasing, it is enough to prove that $\gamma \mapsto h(\gamma)$ is convex. First, let us observe that since $C(\gamma)$ belongs to a fixed circle, the dependence $\gamma \mapsto C(\gamma)$ is continuous. As a direct consequence $\gamma \mapsto h(\gamma)$ is continuous.

Consider $0\leq \gamma_1 < \gamma_2 \leq \pi/6$. Then $h(\gamma_1), h(\gamma_2)$ are the parallel sides of a trapezoid formed by $C(\gamma_1), C(\gamma_2)$ and their projections on $AB$. Since $\gamma_1<\gamma_2$ it follows at once that $AC(\gamma_1)<AC(\gamma_2)$, since by construction $C(\gamma_2)$ has both its coordinates greater in absolute value than those of $C(\gamma_1)$. The line $AC(\frac{\gamma_1+\gamma_2}{2})$ is the bisector of the angle $\angle C(\gamma_1)AC(\gamma_2)$, therefore it intersects $C(\gamma_1)C(\gamma_2)$ in a point $X$ which is closer to $C(\gamma_1)$ than to $C(\gamma_2)$. The point $C(\frac{\gamma_1+\gamma_2}{2})$ is the intersection of $AX$ with the circle $D(1-t)$. Denoting $M$ the midpoint of $C(\gamma_1)C(\gamma_2)$ and $X',M'$ the projections of $X$ and $M$, respectively, on $AB$ we have
\[ h\left(\frac{\gamma_1+\gamma_2}{2}\right) \leq XX' < MM' = \frac{h(\gamma_1)+h(\gamma_2)}{2}.\]
The geometric aspects of the proof are illustrated in Figure \ref{fig:convexity}. 
Since $h$ is also continuous and increasing it follows that $h$ is convex. The previous inequality is strict, as soon as $\gamma_1\neq \gamma_2$, therefore $h$ is strictly convex. In conclusion $\alpha = \arcsin \circ h$ are strictly convex and strictly increasing on $[0,\pi/6]$.

The inequality $\alpha(\gamma)\leq \gamma$ follows directly from $\gamma \in [0,\pi/6]$ and $1-t\geq \sqrt{3}/3$. This implies that $C(\gamma)$ is closer to $A$ than $B$ in the triangle $AYB$, defined previously.
 \hfill $\square$
 
\begin{rem}
	The configuration of triangle $ABC(\gamma)$ in Proposition \ref{prop:convexity} allows the explicit computation of $\alpha(\gamma)$ in terms of $\gamma$ and $t$. Indeed, straightforward computations give $AC(\gamma) = \cos \gamma-\sqrt{(1-t)^2-\sin^2\gamma}$, $h(\gamma) = AC(\gamma)\sin \gamma$ and $\alpha(\gamma) = \arcsin \left( \frac{ h(\gamma)}{1-t}\right)$. 
	\label{rem:explicit-formula}
\end{rem}

\section{Areas and perimeters of inner parallel sets}
\label{sec:area-perim}

The Reuleaux triangle $R$ of width $1$ minimizes the area at fixed constant width $1$. It is natural to conjecture that the same happens for areas of inner parallel sets. Denote by $\Omega_{-t}$ the inner parallel set at distance $t$ from the boundary of $\Omega$, as recalled in  \eqref{eq:def-inner-parallel}. Supposing that $\Omega$ has constant width equal to $1$ here are some well known facts for the extremal cases $t\in \{0,r(R)\}$ where $r(K)$ denotes the inradius of $K$.
\begin{itemize}[noitemsep]
	\item $t=0$: In view of Blaschke-Lebesgue theorem the area of $\Omega_0$ is equal to the area of $\Omega$ and is minimal for the Reuleaux triangle.
	\item $t=r(R)$: It is known that the Reuleaux triangle minimizes the inradius (see \cite{yaglom-boltjanskii}) among shapes with fixed constant width. Therefore the inradius of $\Omega$ is at least equal to $r(R)$, showing that $\Omega_{-r(R)}$ is well defined and non-degenerate. Choosing $t=r(R)$ equal to the inradius of the Reuleaux triangle of width $1$, obviously $|\Omega_{-r(R)}|\geq |R_{-r(R)}| =0$. 
\end{itemize}
In the following, we show that the result holds for all inner parallel sets. First, it is shown that the perimeter of inner parallel sets is minimized by the Reuleaux triangle.

\begin{thm}
	Suppose $\Omega$ has constant width equal to $1$. Given $t\geq 0$, the perimeter of the inner set $|\partial \Omega_{-t}|$ is minimal when $\Omega$ is the Reuleaux triangle $R$ of width $1$. If for some fixed $t>0$ the equality $|\partial \Omega_{-t}|=|\partial R_{-t}|$ holds, then $\Omega=R$.
	\label{thm:perim-parallel}
\end{thm}

\emph{Proof:} Any constant width shape can be approximated arbitrarily well by a Reuleaux polygon (see Proposition \ref{prop:properties}). Therefore, we prove the result for the case of Reuleaux polygons and the general result will follow by a density argument. All constant width shapes are supposed to have width $1$ in the following.

{\bf Step I. Structure of the inner parallel sets for a Reuleaux polygon.} Any Reuleaux polygon has an odd number of sides $n =2k+1\geq 3$. In the following the indices are taken modulo $n$. Following the description in \cite{Firey1960} we represent a Reuleaux polygon $\Omega$ as an intersection of $n$ disks $D_0,...,D_{n-1}$ of radius $1$ centered at points $C_0,...,C_{n-1}$. Moreover, the inner parallel set $\Omega_{-t}=\{ x \in \Omega : d(x,\partial \Omega) \geq t\}$, defined for $t\geq $ smaller than the inradius of $\Omega$, is exactly the intersection of the disks $D_0(1-t),...,D_{n-1}(1-t)$ of radius $1-t$ centered at the same points $C_0,...,C_{n-1}$. We would like to be able to compute or estimate the perimeter of the inner parallel set $\Omega_{-t}$ in terms of the geometry of $\Omega$. At a first sight, we notice that $\partial \Omega_{-t}$ is a union of arcs of circles $\Gamma_i(1-t)$ of radii $(1-t)$ centered at $C_i$. However, the number of arcs and the way their length is computed may change with $t$ as underlined below.

For $t$ small enough, the vertices of $\partial \Omega_{-t}$ lie on the bisector lines corresponding to the vertices of the Reuleaux polygon $\Omega$, since these vertices lie at equal distance from adjacent arcs in the boundary. Therefore, if any two adjacent bisector lines meet before reaching the incenter (the center of the inscribed disk), the corresponding arc vanishes in $\partial \Omega_{-}$.

Let us consider the following iterative process, using the notion of disk-polygons introduced in Section \ref{sec:disk-poly}. At every step $k$ of the iterative process, we will keep track of a set of indices $I_k \subset \{0,...,n-1\}$ corresponding to the disks appearing in the definition of the Reuleaux polygon which will define some disk-polygon.
\begin{itemize}[noitemsep]
	\item {\bf Initialization.} Initialize $\Omega^0 = \Omega$, $t=0$, $k=0$, $I_0=\{0,1,...,n-1\}$.
	\item {\bf Iteration $k$.} Increase $t$ until two adjacent bisectors of $\Omega^k$ meet. If the incenter is reached then stop. Otherwise, one arc $\Gamma_i(1-t)$ of $\Omega^k_{-t}$ (which coincides with $\Omega_{-t}$) is reduced to a point. Remove the corresponding disk $D_i$ from the definition of $\Omega^k$ and define $I_{k+1} = I_k \setminus \{i\}$, $\Omega^{k+1} = \bigcup_{i \in I_{k+1}} D_i$. 
	
	Of course, $\Omega^k_{-t}$ is also a inner parallel set of $\Omega^{k+1}$, since the vanishing arc $\Gamma_i(1-t)$ has no correspondent in the boundary of $\Omega^{k+1}$.
	
	If more than one arc vanishes at a given $t$ then remove all corresponding disks from the definition of $\Omega^k$.

\end{itemize}

\begin{figure}
	\begin{center}
	\includegraphics[height=0.4\textwidth]{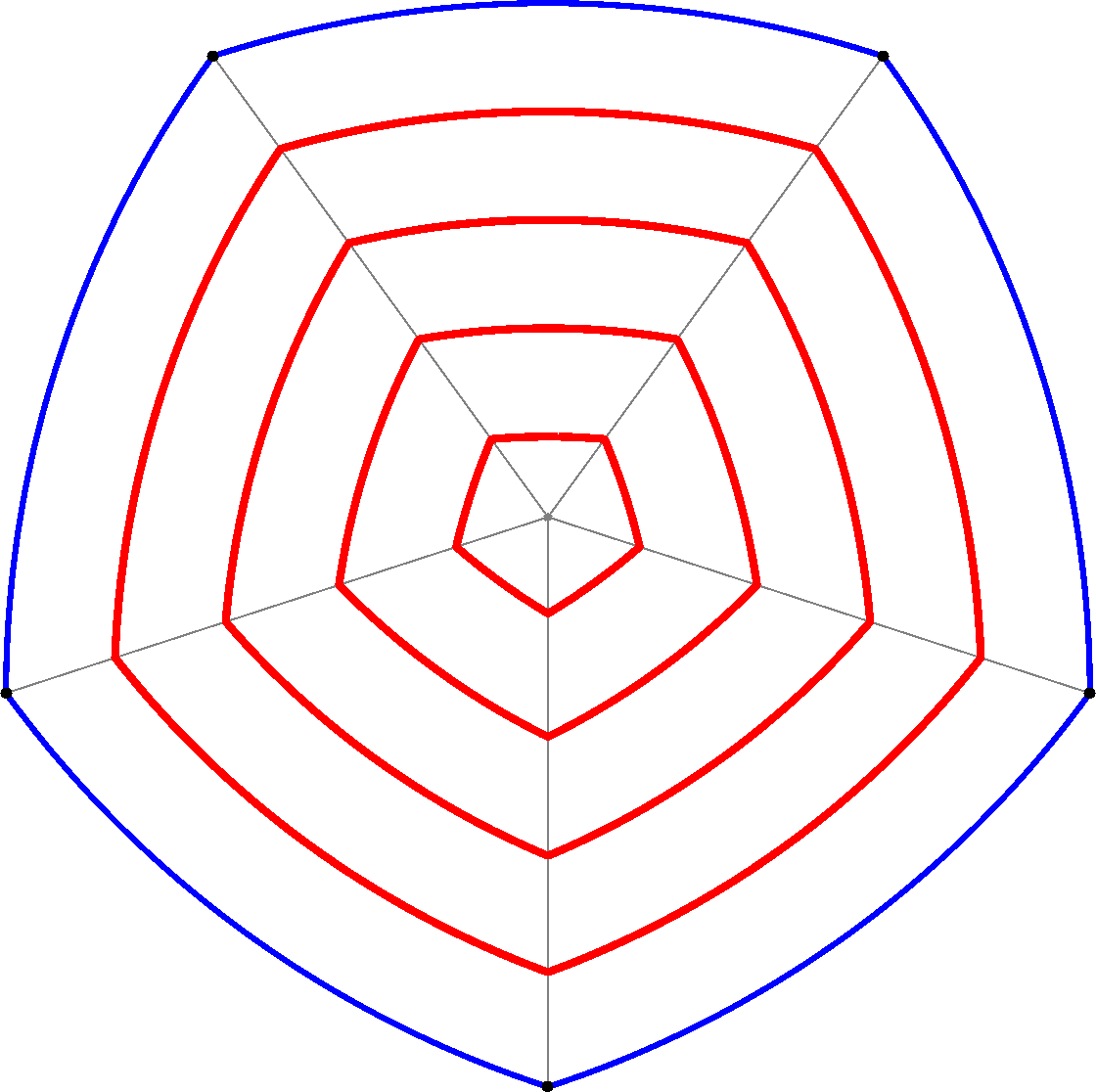}\quad 
	\includegraphics[height=0.4\textwidth]{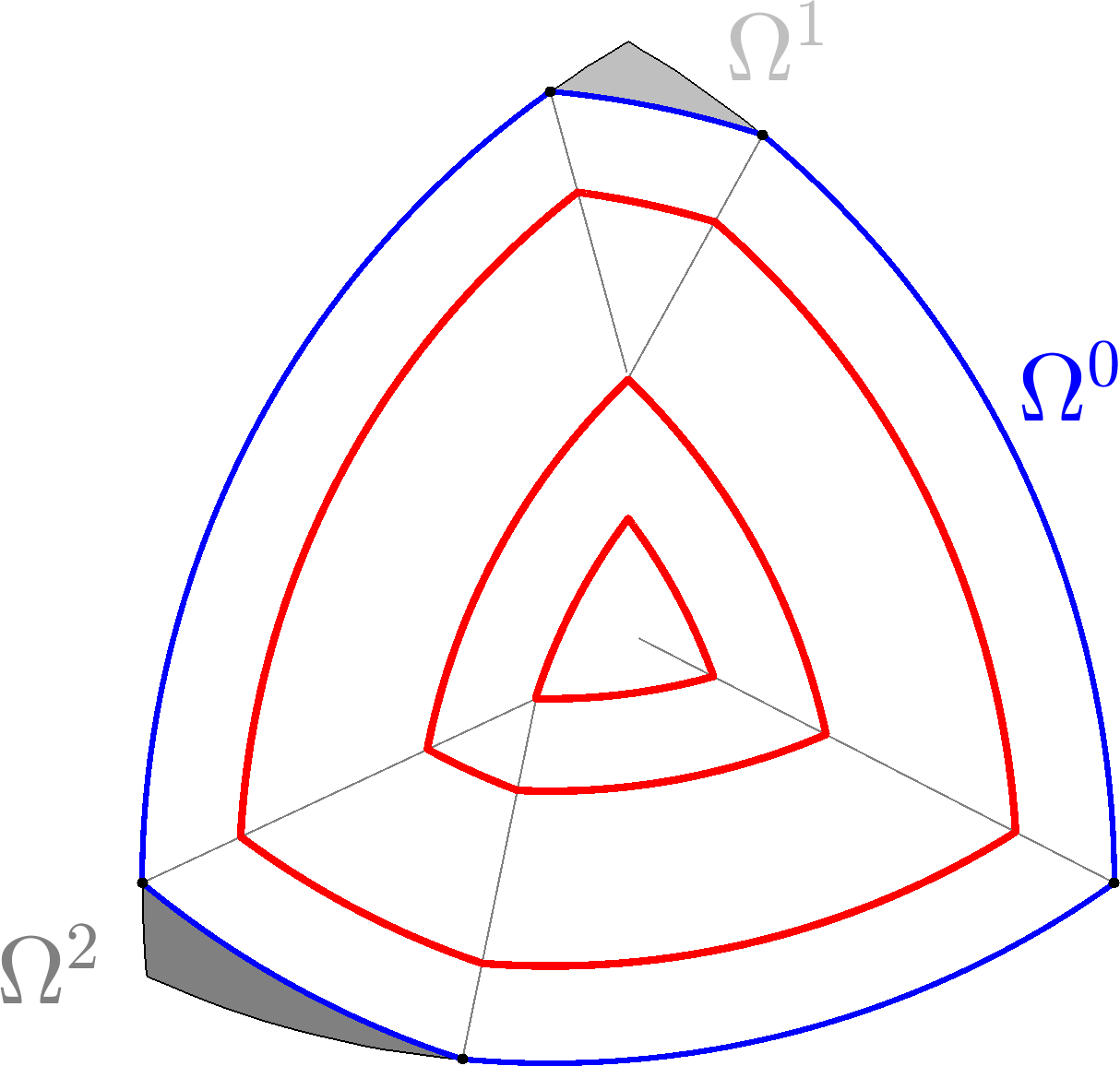}
	\end{center}
	\caption{Examples of inner parallel sets. Left: the regular Reuleaux pentagon has all its inner parallel sets regular. Right: Illustration of the algorithm in \bo{Step I}. Each time two bisectors meet, an edge disappears in the inner parallel set.}
	\label{fig:level-sets-examples}
\end{figure}

The previously defined iterative process ends when the incenter is reached and it has a finite number of steps, since there are finitely many disks involved in the definition of the initial Reuleaux polygon $\Omega$. For example, in the case of a regular Reuleaux polygon, the initialization is enough to reach the incenter, by symmetry. In general situations, multiple iterations may be needed. See Figure \ref{fig:level-sets-examples} for an illustration.

The previous construction is needed because we will be able to compute explicitly the perimeter of the inner parallel set associated to a disk-polygon $\Omega^k$ only in the initial phase, before any two consecutive bisectors of $\Omega^k$ meet. We call these inner parallel sets \emph{regular}, as in the computations related to polygons in \cite{kawohl-lachand-robert}.

At the end of the iterative process, we have an increasing family of disk-polygons (w.r.t. inclusion), all containing the initial Reuleaux polygon $\Omega$, such that all inner parallel level sets $\Omega_{-t}$ are among \emph{regular inner parallel sets} of some disk-polygon $\Omega^k$. More precisely, there exist real numbers $0<t_1<...<t_M$ such that:
\begin{itemize}[noitemsep,topsep=0pt]
	\item $\Omega_{-t}$ is regular for $\Omega_0=\Omega$ for $t \in [0,t_1]$
	\item $\Omega_{-t}$ is regular for $\Omega_1$ for $t \in [t_1,t_2]$
	\item $\Omega_{-t}$ is regular for $\Omega_2$ for $t \in [t_2,t_3]$, etc.
\end{itemize}
The method presented implies that for any $t \in [0,r(\Omega)]$ there exists a disk-polygon $\Omega'\supset \Omega$ such that $\Omega_{-t}=\Omega_{-t}'$ and $\Omega'_{-t}$ is a regular inner parallel set in $\Omega'$.

{\bf Step II. Computing the perimeter of a regular inner parallel set of a disk-polygon.} Consider a disk-polygon $\Omega = \bigcap_{i=0}^{N-1} D_i$, where $D_i$ are disks of radius $1$, whose boundary is made of arcs having measures $\theta_i$ and turning angles $\beta_i$, $i=0,...,N-1$. Following the description in Section \ref{sec:disk-poly} and Proposition \ref{prop:disk-poly} we have $\sum_{i=0}^{N-1} (\theta_i+\beta_i)=2\pi$. 

Consider $t^*>0$ such that for $t\in [0,t^*]$ no two consecutive angle bisectors of $\Omega$ meet. Therefore, for $t \in [0,t^*]$ the inner parallel set $\Omega_{-t}$ of the disk-polygon $\Omega$ is regular. 

Suppose that $\Omega$ is one of the disk-polygons $\Omega^k$ derived in {\bf Step I}. Then $\Omega^k$ contains the constant width shape $\Omega$, the initial Reuleaux polygon. Thus, following Proposition \ref{prop:disk-poly} the perimeter of $\Omega^k$ is larger than $\pi$ and we have $\sum_{i=0}^{N-1} \theta_i \geq \pi$, implying that $\sum_{i=0}^{N-1} \beta_i \leq \pi$. Moreover, since $\Omega^k$ has minimal width at least equal to $1$, we have $\beta_i \in [0,\pi/3]$ for all $i=0,...,N-1$.

\begin{figure}
	\centering 
	\includegraphics[height=0.7\textwidth]{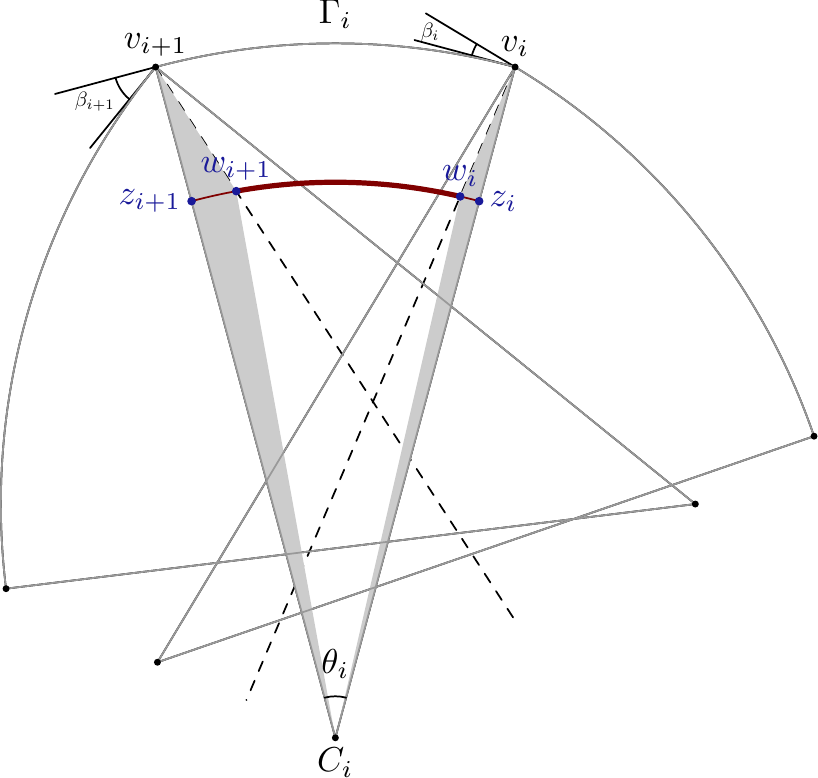}
	\caption{Computation of the perimeter of a regular inner parallel set for a disk-polygon. The shaded triangles have the same configuration as $\Delta ABC(\gamma)$ in Proposition \ref{prop:convexity} and Figure \ref{fig:convexity}.}
	\label{fig:disk-poly}
\end{figure}

Given an arc $\Gamma_i$ on the boundary of $\Omega$ we compute its contribution to the boundary of $\Omega_{-t}$. Consider $S_i$ the circular sector having angle $\theta_i$ and center $C_i$ associated to $\Gamma_i$ determined by $\Gamma_i$ and the normals and its endpoints $v_i,v_{i+1}$. The angle bisectors of $\Omega$ at $v_i$ and $v_{i+1}$ are contained in $S_i$ and make angles $\beta_i/2,\beta_{i+1}/2$, respectively, with the rays of the sector. See Figure \ref{fig:disk-poly}.

Given $t \in [0,t^*]$, denote by $w_i,w_{i+1}$ the intersection of the bisectors at vertices $v_i, v_{i+1}$ with the circle of center $C_i$ and radius $1-t$. Denote by $z_i, z_{i+1}$ the intersections of the same circle with the boundary rays of the sector $S_i$. See Figure \ref{fig:disk-poly} for an illustration. Note that triangles shaded triangles in Figure \ref{fig:disk-poly} have the same configuration as in Figure \ref{fig:convexity} and Proposition \ref{prop:convexity}. More precisely $\Delta v_iC_iw_i\equiv \Delta ABC(\beta_i/2)$ and $\Delta v_{i+1}C_iw_{i+1}\equiv \Delta ABC(\beta_{i+1}/2)$. Then the contribution of $\Gamma_i$ to the perimeter of $\Omega_{-t}$ is equal to 
\[ |\arc{w_iw_{i+1}}|=|\arc{z_iz_{i+1}}|-|\arc{w_iz_i}|-|\arc{w_{i+1}z_{i+1}}|=(1-t)(\theta_i-\alpha_t(\beta_i/2)-\alpha_t(\beta_{i+1}/2)),\]
where the function $\alpha_t$ is the one defined in Proposition \ref{prop:convexity} and explicited in Remark \ref{rem:explicit-formula}. 

Thus, we obtain that for all $t \in [0,t^*]$ we have
\begin{equation} |\partial \Omega_{-t}| = (1-t)\left(\sum_{i=0}^{N-1} \theta_i - 2 \sum_{i=0}^{N-1} \alpha_t(\beta_i/2)\right) .
\label{eq:regular-perim}
\end{equation}
Note that the function $\alpha_t$ defined in Proposition \ref{prop:convexity} is convex on $[0,\pi/6]$, increasing and depends implicitly on $t$. Moreover, we work with $t\leq r(R)=\frac{3-\sqrt{3}}{3}$ (the inradius of the Reuleaux triangle of width $1$).

In the case of regular Reuleaux polygons, in particular for the Reuleaux triangle, all inner parallel sets are regular, since the angle bisectors meet at the incenter, because of the symmetry. Therefore, if $\Omega$ is a regular Reuleaux polygon, for all $t \leq r(\Omega)$ we have
\begin{equation} |\partial \Omega_{-t}| = (1-t)\left(\pi - 2n \alpha_t\left(\frac{\pi}{2n}\right)\right),
\label{eq:regular-reuleaux-perim}
\end{equation}
where we used the fact that $\theta_i=\beta_i = \pi/n$.

{\bf Step III. Optimality of the Reuleaux triangle via convexity arguments.} 

For $t\in [0,t^*]$ we define the convex and strictly increasing function (see Proposition \ref{prop:convexity}) 
\begin{equation} f:[0,\pi/3] \to \Bbb{R}, f(\beta) = 2\alpha_t(\beta/2).
\label{eq:alpha-f}
\end{equation}
For fixed $t$ and a fixed disk-polygon defined as in the previous step, there exists $\varepsilon\in [0,\pi]$ such that $\sum_{i=0}^{N-1} \theta_i = \pi+\varepsilon$ and $\sum_{i=0}^{N-1} \beta_i = \pi-\varepsilon$. Thus, \eqref{eq:regular-perim} becomes
\begin{equation}
|\partial \Omega_{-t}| = (1-t)(\pi+\varepsilon) - (1-t) \sum_{i=0}^{n-1} f(\beta_i).
\label{eq:perim-lower-bound}
\end{equation}


Consider the maximization problem 
\begin{equation}
  \begin{array}{c}
  M(\varepsilon) = \max \left[f(\beta_0)+ f(\beta_1)+...+f(\beta_{N-1})\right] \\
  \text{ such that } \beta_i \in [0,\pi/3], i=0,...,N-1, \beta_0+...+\beta_{N-1}=\pi-\varepsilon.
  \end{array}
  \label{eq:max-sum-f}
\end{equation}
Since $f$ is convex and $C^1$ on $[0,\pi/3]$, if $a<b\in (0,\pi/3)$, the function $g(t) = f(a-t)+f(b+t)$ is strictly increasing in a neighborhood of $0$. Indeed, we have
\[ g'(t)=f'(b+t)-f'(a-t)>0,\]
since $f'$ is strictly increasing. Therefore if $t>0$ is small enough, replacing $(a,b)$ with $(a-t,b+t)$ preserves the constraints and increases the value of the objective function. Furthermore, we can choose $t$ such that either $a-t=0$ or $b+t=\pi/3$. Therefore, a maximizer for \eqref{eq:max-sum-f}, which exists by classical compactness arguments, has at most one $\beta_i \in (0,\pi/3)$ and all other $b_j$ for $j\neq i$ belong to $\{0,\pi/3\}$. Moreover, the maximal value $M(\varepsilon)$ is clearly strictly decreasing in $\varepsilon$ since $f$ defined in \eqref{eq:alpha-f} is strictly increasing. Thus $M(\varepsilon)\leq M(0)=3f(\pi/3)$.

Therefore, \eqref{eq:perim-lower-bound} and the previous considerations related to problem \eqref{eq:max-sum-f} imply
\begin{equation} |\partial \Omega_{-t}| \geq (1-t)(\pi+\varepsilon)-(1-t)M(\varepsilon)  \geq (1-t)(\pi-3f(\pi/3)) = |\partial R_{-t}|.
\label{eq:final-estimate}
\end{equation}
Thus, for the given $t$ the Reuleaux triangle has a smaller perimeter than $|\partial \Omega_t|$.

To conclude, given $t\leq r(R)$ any inner parallel set of a Reuleaux polygon is a \emph{regular} inner parallel set of a disk-polygon containing it, described in {\bf Step I}. Computing its perimeter as shown in \bo{Step II} and using the estimate in \bo{Step III} shows that $|\partial \Omega_{-t}|\geq |\partial R_{-t}|$. Using the density of Reuleaux polygons in the class of constant width shapes (see Proposition \ref{prop:properties}) gives the same estimate for a general shape having constant width $1$. 

If for some $t>0$ we have $|\partial\Omega_{-t}|=|\partial R_{-t}|$, then equality holds in \eqref{eq:final-estimate}. Therefore $\varepsilon=0$ and $(\beta_i)_{i=0}^{n-1}$ must solve \eqref{eq:max-sum-f}, which shows that the turning angles $\beta_i$ correspond to a Reuleaux triangle.
\hfill $\square$

Now we are ready to prove a similar result for the areas of the inner parallel sets.

\begin{thm}
	Suppose $\Omega$ has constant width equal to $1$. Given $t\in [0,r(R)]$, the area of the inner set $|\Omega_{-t}|$ is minimal when $\Omega$ is the Reuleaux triangle $R$ of width $1$. If for some fixed $t\geq 0$ the equality $|\Omega_{-t}|=|R_{-t}|$ holds, then $\Omega=R$.
	\label{thm:area-parallel}
\end{thm}

\emph{Proof:} We suppose that $\Omega$ is a Reuleaux polygon whose edges are arcs of circles of radius $1$ centered in $C_i$ and having arc measures $\theta_i$, $i=1,...,k$. Of course, we have the well known property $\sum_{i=1}^k \theta_i=\pi$ (see Proposition \ref{prop:properties}).  

Denote by $A_\Omega(t) = |\Omega_{-t}|$ the area and $P_\Omega(t) = |\partial \Omega_{-t}|$ the perimeter of the inner parallel set at distance $t$. One can see that changing $t$ induces a normal movement of the boundary with uniform speed. Recall that the shape derivative of the area is given by $|\omega|'(V)= \int_{\partial \omega} V\cdot n$, where $n$ is the outer normal. This is classical and can be found, for example, in \cite[Section 2.5]{sokolowski-zolesio}, \cite[Chapter 5]{henrot-pierre-english}. In the case of inner parallel sets of Reuleaux polygons, direct computations can be made like in \cite{Firey1960}. Therefore, classical shape derivative formulas imply that $A_\Omega'(t) = -P_\Omega(t)$.

Let $r=r(R)$ be the inradius of the Reuleaux triangle of width $1$, which is minimal among shapes of constant width. Thus, for $t \geq r$ the inner parallel set $\Omega_{-t}$ is non-trivial. Moreover, from Theorem \ref{thm:perim-parallel} we have $P_\Omega(t)\geq P_R(t)$, which gives
\[ A_\Omega'(t) = -P_\Omega(t)  \leq - P_R(t) = A_R'(t).\]
Therefore $t\mapsto A_\Omega(t)-A_R(t)$ is decreasing. Comparing values in $t\leq  r$ we obtain
\begin{equation} A_\Omega(t)-A_R(t) \geq A_\Omega(r)-A_R(r) = A_\Omega(r)\geq  0.
\label{eq:ineq-area}
\end{equation}
Using the density of Reuleaux polygons in the class of constant width shapes (see Proposition \ref{prop:properties}) finishes the proof.

If for some $t\geq 0$ we have $A_\Omega(t) = A_R(t)$ then \eqref{eq:ineq-area} implies that $A_\Omega(r)=0$, i.e. the inradius of $\Omega$ smaller or equal to the inradius of the Reuleaux triangle. This implies $\Omega=R$.
\hfill $\square$

\begin{cor}
	(Blaschke-Lebesgue theorem) Applying Theorem \ref{thm:area-parallel} for $t=0$ shows that the Reuleaux triangle minimizes the area among shapes of constant width.
	\label{cor:bl}
\end{cor}

\emph{Proof:} The result is obvious from Theorem \ref{thm:area-parallel}. Nevertheless, we must underline that the Blaschke-Lebesgue theorem was not used in any of the preceding results. We only used the minimality of the inradius for the Reuleaux triangle, which is proved directly, via classical geometry arguments in \cite[Exercise 7-14]{yaglom-boltjanskii}, for example. \hfill $\square$

\begin{rem}
	Corollary \ref{cor:bl} shows that, in a certain sense, the minimality of the inradius for the Reuleaux triangle among shapes of constant width is a stronger result which implies the Blaschke-Lebesgue theorem, when combined with the previous independent results shown in the paper.
\end{rem}

\begin{rem}
	Theorem \ref{thm:area-parallel} can also be deduced from the results of \cite{disk-polygons-BL} where it is shown that if $\Omega$ is a disk-polygon made with disks of unit radius whose centers are at distance at most $d\in [1,\sqrt{3}]$ apart, then $|\Omega| \geq \Delta(d)$, where $\Delta(d)$ is the disk polygon obtained putting the centers at the vertices of an equilateral triangle of side $d$. 
	
	Rescaling we find that the vertices of a disk polygon with disks having rays $r \in [\frac{\sqrt{3}}{3},1]$ whose centers are at distance at most $1$ apart, the same result holds, comparing the inner parallel areas of a Reuleaux polygon and the Reuleaux triangle.
\end{rem}

\section{Blaschke-Lebesgue Theorem for Cheeger sets}
\label{sec:cheeger}

For a bounded and convex domain $\Omega$, the associated Cheeger constant is defined by 
\begin{equation}\label{eq:cheeger}
h(\Omega)  = \min_{E\subset \Omega} \frac{|\partial E|}{|E|},
\end{equation}
where $|\partial E|$ denotes the perimeter of the set $E$, which may be assumed convex when $\Omega$ is convex.
This notion was introduced by Cheeger in \cite{cheeger} in order to give a geometric lower bound for the first eigenvalue of the Dirichlet Laplacian. The Cheeger constant, although defined using geometric quantities in \eqref{eq:cheeger}, can be interpreted as the first eigenvalue of the $1$-Laplacian \cite{kawohl-fridman}.

For convex planar domains Lachand-Robert and Kawohl show in \cite{kawohl-lachand-robert} that $h(\Omega)$ can be characterized using areas of parallel inner sets. Indeed, we have
\begin{equation}
 h(\Omega) = 1/t, \text{ where } |\Omega_{-t}| = \pi t^2.
 \label{eq:cheeger-char}
\end{equation}

Recently, in \cite{BHL18} the minimization of the eigenvalues of the Dirichlet-Laplace operator was studied for domains having a diameter constraint. The optimal shapes are shapes of constant width. The maximization of the eigenvalues of the Dirichlet-Laplace operator also makes sense in the class of shapes of constant width (see \cite[Section 4.2]{Bogosel-Convex}). Numerical simulations shown in \cite[Section 4.2]{Bogosel-Convex} indicate that the Reuleaux triangle is likely to be the constant width shape which maximizes these eigenvalues. For now, no proof of this fact is known, as recalled in \cite[Section 5]{henrot-lucardesi}. This fact is further motivated by the fact that the first eigenvalue of the $\infty$-Laplacian, which is the inverse of the inradius, is also maximized by the Reuleaux triangle.

It is natural, therefore, to conjecture that the Cheeger constant is maximized by the Reuleaux triangle, in the class of shapes of constant width. This result was proved recently in \cite{henrot-lucardesi} using techniques from shape optimization, notably the optimality condition verified by a minimizer. 

The proof given below is quite straightforward, following the characterization \eqref{eq:cheeger-char} and the result of Theorem \ref{thm:area-parallel}. The existence of constant width shapes maximizing the Cheeger constant is straightforward due to classical compactness arguments among convex sets and the continuity of the Cheeger constant. The proof given below is, however, based only on a direct comparison with the Reuleaux triangle.

\begin{thm}
	(Blaschke-Lebesgue Theorem for the Cheeger constant) The Reuleaux triangle $R$ maximizes $h(\Omega)$ when $\Omega$ has fixed constant width. Moreover, $R$ is the unique maximizer.
	\label{thm:cheeger}
\end{thm}

\emph{Proof:} Let $h(\Omega)=1/t$ be the Cheeger constant of the constant width shape $\Omega$. Then $|\Omega_{-t}|=\pi t^2$, in view of \eqref{eq:cheeger-char}. If $R$ denotes the Reuleaux triangle of the same width, then Theorem \ref{thm:area-parallel} shows that $|R_{-t}|\leq |\Omega_{-t}|=\pi t^2$.

Consider now the Cheeger constant $h(R) = \frac{1}{t^*}$ of the Reuleaux triangle $R$. Then \eqref{eq:cheeger-char} gives $|R_{-t^*}|=\pi (t^*)^2$. 

\begin{figure}
	\centering
	\includegraphics[height=0.35\textwidth]{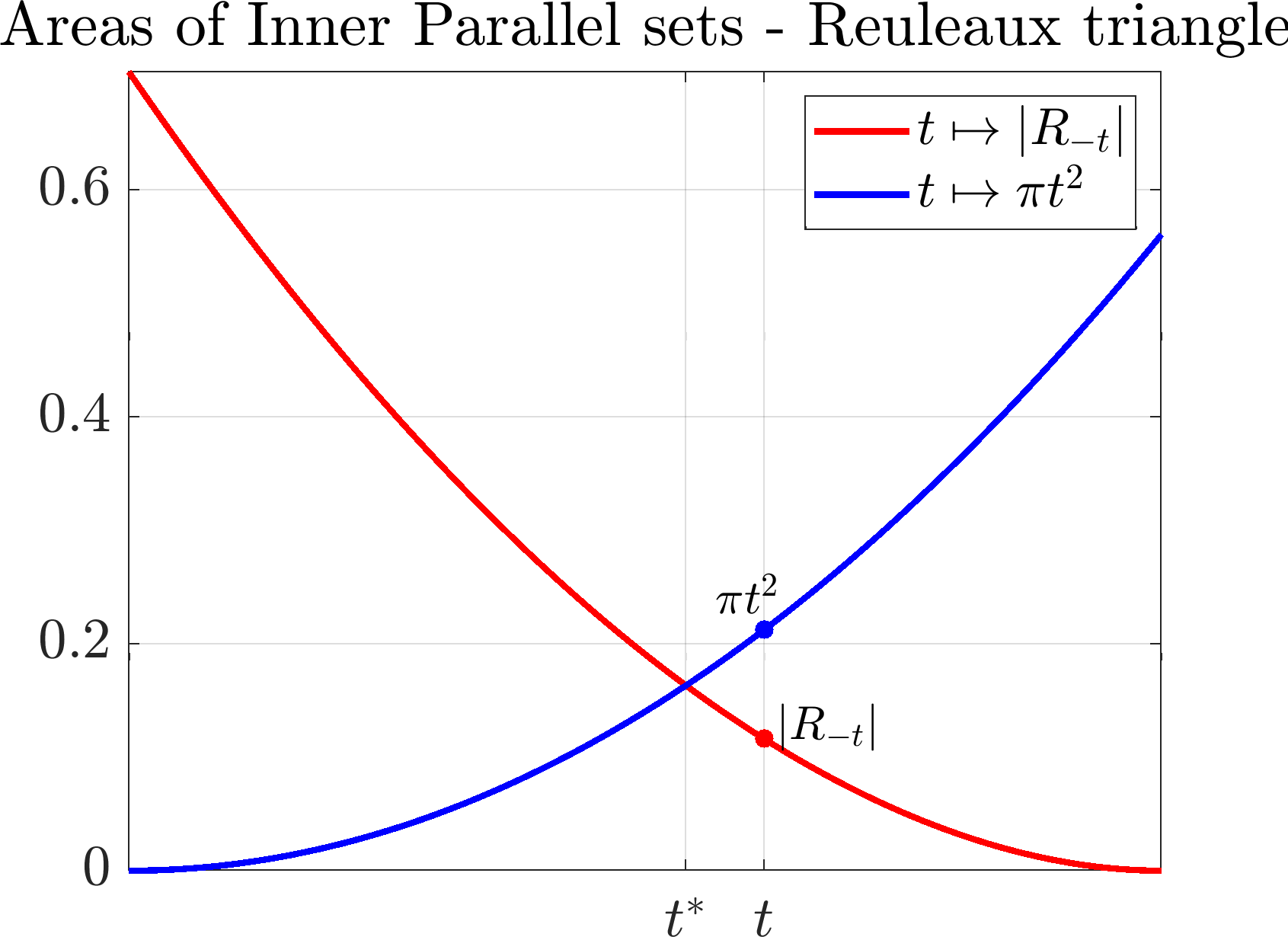}\quad
	\includegraphics[height=0.35\textwidth]{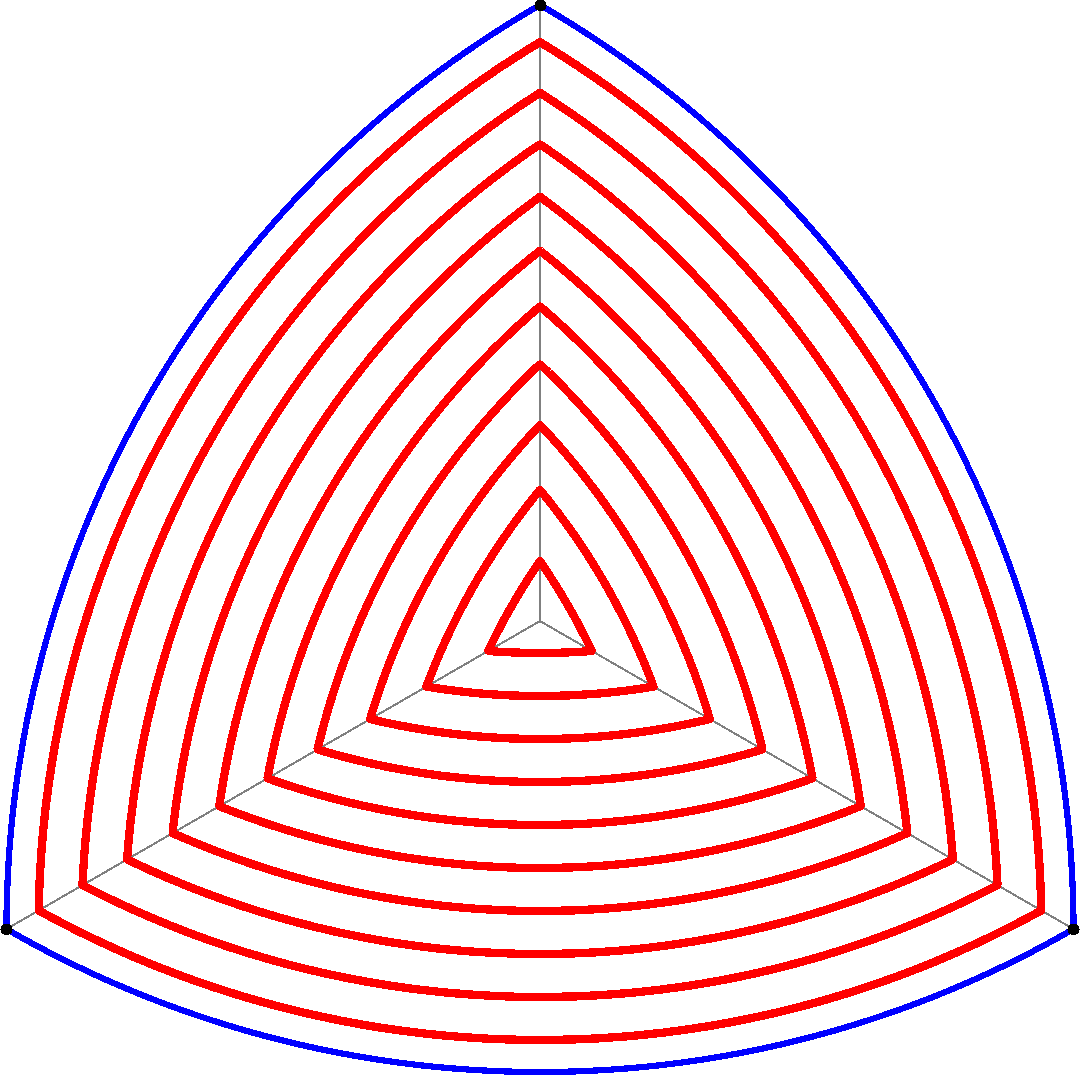}
	\caption{Left: graphical representation for areas of inner parallel sets of the unit Reuleaux triangle $R$ versus the graph of $t \mapsto \pi t^2$. The $x$-coordinate $t^*$ of the point of intersection is the inverse of the Cheeger constant $h(R)$. The graph gives a geometrical proof of the implication $\pi t^2 \geq |R_{-t}| \Longrightarrow t\geq t^*$. Right: graphical representation of a few parallel inner sets of $R$.}
	\label{fig:ReuleauxCheeger}
\end{figure}

Recall that the areas of parallel sets $t\mapsto |R_{-t}|$ are strictly decreasing for $t\in [0,r(R)]$ while $t\mapsto \pi t^2$ is strictly increasing. The value of $t$ for which the two functions are equal is precisely $t^*$. However, since $|R_{-t}|\leq \pi t^2$, we find that $t$ must lie in the region where $|R_{-t}|$ is below $\pi t^2$, i.e. $t\geq t^*$. See the illustration in Figure \ref{fig:ReuleauxCheeger}. This shows that
\[ h(R) = \frac{1}{t^*} \geq \frac{1}{t} = h(\Omega).\]

If $h(R)=h(\Omega)$ then, following the inequalities shown above, we find that $|R_{-t}|=|\Omega_{-t}|$, which according to Theorem \ref{thm:area-parallel} implies that $\Omega=R$.
\hfill $\square$

{\bf Acknowledgments.} The author thanks Dorin Bucur for suggesting the question solved in Theorem \ref{thm:cheeger} back in 2016. This work was supported by the ANR Shapo program (ANR-18-CE40-0013).

\bibliographystyle{abbrv}
\bibliography{./biblio}

\bigskip
\small\noindent
Beniamin \textsc{Bogosel}: Centre de Math\'ematiques Appliqu\'ees, CNRS,\\
\'Ecole polytechnique, Institut Polytechnique de Paris,\\
91120 Palaiseau, France \\
{\tt beniamin.bogosel@polytechnique.edu}\\
{\tt \nolinkurl{http://www.cmap.polytechnique.fr/~beniamin.bogosel/}}

\end{document}